\documentclass[reqno,12pt]{amsart}
\usepackage{amsmath,amsthm,amssymb}
\usepackage{amsfonts}
\usepackage{hyperref}
\usepackage{bbm}

\date{}

\newtheorem{theorem}{Theorem}[section]

\newtheorem{lemma}[theorem]{Lemma}

\newtheorem{proposition}[theorem]{Proposition}

\newtheorem{remark}[theorem]{Remark}
\numberwithin{equation}{section}

\begin{document}

\centerline{\large {\bf
 Sobolev--Kantorovich inequalities under $CD(0, \infty)$ condition
}}

\vskip .3in

\centerline{Vladimir I. Bogachev}
\vskip .1in

\centerline{{\it Department of Mechanics and Mathematics, Moscow State University}}
\centerline{{\it 119991 Moscow, Russia}}
\centerline{{\it National Research University Higher School of Economics}}
\centerline{{\it Faculty of Mathematics, Usacheva 6}}
\centerline{{\it 119048 Moscow, Russia}}
\centerline{{\it  vibogach@mail.ru}}

\vskip .2in

\centerline{Alexander V. Shaposhnikov}
\vskip .1in

\centerline{{\it Faculty of Mathematics, University of Bielefeld}}
\centerline{{\it  D-33615 Bielefeld, Germany}}
\centerline{{\it  shal1t7@mail.ru, ashaposh@math.uni-bielefeld.de}}

\vskip .2in

\centerline{Feng-Yu Wang}
\vskip .1in

\centerline{{\it Center for Applied Mathematics,  Tianjin University}}
\centerline{{\it Tianjin 300072, China}}
\centerline{{\it Department of Mathematics, Swansea University}}
\centerline{{\it Bay Campus, SA1 8EN, UK}}
\centerline{{\it wangfy@tju.edu.cn}}

\vskip .2in
\centerline{{\bf Abstract}}

We refine and generalize several interpolation inequalities
bounding the $L^p$ norm of a probability density with respect to the reference measure
$\mu$ by its Sobolev norm and the Kantorovich distance to $\mu$
on a smooth weighted Riemannian manifold satisfying $CD(0, \infty)$ condition.

\vskip .1in



\vskip .2in

\section{Introduction}

In the last decade there has been an increasing interest in functional inequalities
relating Sobolev norms  with  certain other norms and quantities such as
entropy and Kantorovich distances, e.g. see \cite{BakryGL}, \cite{WFBook}, \cite{W13}.
In this paper we discuss some interpolation inequalities which can be viewed as analogs of the classic Hardy--Landau--Littlewood inequality
$$
\|f'\|_{L^1}^2 \le C\|f\|_{L^1}\|f''\|_{L^1}
$$
as well as the celebrated Otto--Villani HWI inequality
$$
\int f\log f \,d\mu \leq \biggl(\int \frac{|\nabla f|^2}{f}\,d\mu \biggr)^{1/2}
W_{2}(\mu, f\mu).
$$
It was noticed in \cite{BS15}, \cite{BWS15} that the Hardy--Landau--Littlewood inequality can be written  as
$$
\|f\|^{2}_{L^1}\le C\|\nabla f\|_{L^1}\, \|f\|_{K},
$$
where $f$ belongs to the usual Sobolev class
$W^{1,1}$ and has zero integral and $\|f\|_K$ denotes the  Kantorovich norm
of the signed measure $f\,dx$ with zero value on the whole space.
Recall that the Kantorovich norm of a signed measure $\mu$ on $\mathbb{R}^n$ with $\mu(\mathbb{R}^n)=0$ integrating Lipschitz
functions is defined by
$$
\|\mu\|_K=\sup\biggl\{ \int f\, d\mu\colon f\in C_b^\infty(\mathbb{R}^n), \ |\nabla f|\le 1\biggr\}.
$$
In this form this inequality admits natural
multidimensional extensions such as the bound
$$
\|\mu\|^2 \le C \|D\mu\|\, \|\mu\|_K
$$
established in \cite[Theorem 1]{BS15} for signed Borel measures $\mu$ on $\mathbb{R}^n$ with
$\mu(\mathbb{R}^n)=0$ possessing a density of class BV, where $\|\mu\|$ is the total variation of $\mu$ and
$\|D\mu\|$ is the total variation of the vector measure $D\mu$ that is the distributional derivative of~$\mu$.
Next, in our note \cite{BWS15} a dimension-free
version of this bound employing probability reference measures was established, in particular, for the standard Gaussian
measure $\gamma_d$ on $\mathbb{R}^d$ we proved that
$$
\|f\|_{L^1(\gamma_d)}^2\le 2 \|\nabla f\|_{L^1(\gamma_d)} \|f\gamma_d\|_K
$$
for all smooth functions $f\in L^1(\gamma_d)$ with zero integral against~$\gamma_d$.

There are also extensions closely related to the inequalities
arising in the study of some evolution equations (e.g.,  the Cahn--Hilliard model) established
by  Cinti, Kohn, and Otto \cite{CintiOtto},  \cite{KohnOtto}.
In particular, it was proved in  \cite[Proposition 1.3]{CintiOtto}  that for any periodic smooth probability density $f$ on~$[0, 1]^n$
one has
$$
\|(f - C)_{+}\|^{\theta}_{r}
\leq C \|\nabla f\|_{1}W_{2}(dx, f\,dx),
$$
where
$$
r = \frac{3n  + 2}{3n},\ \theta = \frac{3n + 2}{2n}
$$
and
$W_{2}(dx, f\,dx)$ denotes the Kantorovich distance of order $2$ between
the probability measures $dx$ and $f\,dx$ on~$[0, 1]^n$.
This result was generalized by Ledoux \cite{Led}
in the setting of non-negatively curved weighted manifolds.
Let $(M, g)$ be a complete connected $n$-dimensional Riemannian
manifold with the Riemannian volume $dx$ and let
$\mu$ be a probability measure on $M$ with a smooth density
 with respect to $dx$. The term ``smooth function'' will mean below
 a function of class $C_b^\infty(M)$.
The symbol $\|\,\cdot\,\|_q$ will refer to the norm in $L^q(\mu)$.
 In the ``finite-dimensional'' case
Ledoux established the following theorem.

\begin{theorem} {\rm(Theorem 1.1 from  \cite{Led})}
Suppose that $\mu$ satisfies the curvature-dimension condition $CD(0, N)$
for some $N \ge 1$. Given $p, q \geq 1$, there is
a constant $C > 0$ depending only on $p, q, N$
such that for any probability measure
$\nu = f\cdot \mu$ with a smooth density $f$ one has
$$
\|(f - C)_{+}\|^{\theta}_{r}
\leq C \|\nabla f\|_{q}W_{p}(\mu, \nu),
$$
where
$$
r = \frac{1 + \frac{1}{p} + \frac{1}{N}}
{\frac{1}{p} + \frac{1}{q}},
\quad
\theta = r\Bigl(\frac{1}{p} + \frac{1}{q}\Bigr) =
1 + \frac{1}{p} + \frac{1}{N}.
$$
\end{theorem}

In the ``infinite-dimensional'' case $CD(0, \infty)$
the main result from \cite{Led} is summarized in the next theorem.

\begin{theorem}\label{theorem:ledoux_infinite_dimensional}
 {\rm(Theorem 4.1 from  \cite{Led})}
Suppose that $\mu$ satisfies the curvature-dimension condition $CD(0, \infty)$.
Given $1 < q \leq 2$, there is a constant $C > 0$
 depending only on $q$ such that for any
probability measure $\nu = f\cdot \mu$ with a smooth density $f$ one has
$$
\|(f - C)_{+}\|^{3/2}_{r} \leq C\|\nabla f\|_q W_{2}(\mu, \nu),
\quad
r = \frac{3q}{q + 2}.
$$
\end{theorem}

Of course, these bounds extend to densities $f$ from the Sobolev class
$W^{q,1}(\mu)$ defined as the completion of $C_0^\infty(M)$
with respect to the Sobolev norm
$$
\|f\|_{q,1}:=\|f\|_{L^q(\mu)}+ \bigl\|\, |\nabla f|\,\bigr\|_{L^q(\mu)},
$$
associated with $\mu$. Using uniformly Lipschitz bump functions on~$M$
(cf. \cite[Chapter~2]{Aubin}),
it is readily verified that $W^{q,1}(\mu)$ coincides with the class of locally
Sobolev functions on $M$ with finite norm~$\|f\|_{q,1}$.

For the Kantorovich distance of order~$1$ (Kantorovich norm)
 it was proved in \cite{BWS15} that on a smooth weighted
Riemannian manifold $(M, g, \mu)$
satisfying the curvature-dimension condition
$CD(\kappa, \infty)$ with $\kappa \geq 0$
for any smooth function $f$ with zero integral
the following inequality holds:
$$
\|f\|_{1} \leq \inf_{\tau > 0}\Biggl[
\|f\|_{K} \sqrt{\frac{\kappa}{e^{2\kappa \tau} - 1}} +
\|\nabla f\|_{1}\int_{[0, \tau]}
\sqrt{\frac{\kappa}{e^{2\kappa t} - 1}} \,dt
\Biggr],
$$
where for $\kappa = 0, t > 0$ we set
$$
\frac{\kappa}{e^{2\kappa t} - 1} := \frac{1}{2t}.
$$
This bound implies that under the condition $CD(0, \infty)$
 for any two smooth probability densities $f_1, f_2$ the following inequality holds:
$$
\|f_1 - f_2\|^{2}_{1}  \leq 2 \|\nabla f_1 - \nabla f_2\|_1 \cdot
\|f_{1}\cdot \mu - f_{2}\cdot \mu\|_{K}.
$$

In the present paper Theorem \ref{theorem:ledoux_infinite_dimensional} (Theorem 4.1 from  \cite{Led}) is improved in several directions.
We show that a stronger   inequality with an extra
logarithmic factor on the left-hand side holds true for all $q >1$.
Our approach refines the arguments from
\cite{CintiOtto}, \cite{Led}  and combines them with some ideas from our
 short notes \cite{BWS15}, \cite{BWS16}, where a preliminary version of this result was announced without a proof and in smaller generality.
The main result is contained in the next theorem.

\begin{theorem}\label{theorem:main}
Let $(M, g, \mu)$ be a complete connected Riemannian
manifold with a probability measure $\mu$ with a smooth density
satisfying the   condition $CD(0, \infty)$.
For every $q > 1$, there is
a constant $C > 0$ depending only on $q$
such that for any probability measure
$\nu = f\cdot \mu$ with a smooth density~$f$ {\rm(}or $f\in W^{q,1}(\mu)${\rm)}
 one has
$$
\bigl \|(f - C)_{+}\bigl(1 + \log^{\alpha}(1 + (f - C)_{+})\bigr)\bigr\|^{3/2}_{r}
\leq C \|\nabla f\|_{q}W_{2}(\mu, \nu),
$$
where
$$
r = \frac{3q}{q + 2}, \
\alpha = \frac{1}{3}.
$$
The assertion does not hold if   $\log^\alpha(1+\cdot)$ is replaced by any increasing positive function
$\Phi\colon [0,\infty)\to [0,\infty)$ with
$$\lim\limits_{u\to\infty }\Phi(u)\log^{-2/3}(1+u)=\infty.$$
\end{theorem}

\begin{remark}
\rm
It would be interesting to investigate the sharpness of this inequality,
i.e., determine the optimal value of~$\alpha$.
Theorem~\ref{theorem:main} essentially
states that the optimal value of $\alpha$  belongs to
the interval~$[1/3, 2/3]$.
\end{remark}

\begin{remark}
\rm
One of the consequences of Theorem~\ref{theorem:main} is the following
Nash-type inequality for Sobolev functions $f \in W^{2, 1}(\mu)$
which can be of independent interest:
\begin{equation*}
\bigl\| f \log^{1/3}(e + |f|) \bigr\|_{3/2}
\leq C\bigl[\|f\|^{2}_{1} \log(e + \|f\|_1) + \|\nabla f\|^{2}_{2}D(M)^{2}\bigr]^{1/3}\|f\|_{1}^{1/3},
\end{equation*}
where $D(M)$ is the diameter of $M$.
\end{remark}

The paper is organized as follows. In Section 2
we introduce our framework and some basic tools, such as the pseudo-Poincar\'e inequality,
the Kantorovich duality, and the infinite-dimensional Harnack inequality.
In Section~3 we prove Theorem \ref{theorem:main}.
In Section~4 the case of the Kantorovich distance $W_1$ of order~$1$
is discussed. In Section~5 we present some extensions to the case of negatively curved weighted Riemannian manifolds.

\section{Framework and geometric tools}

Throughout this paper, we assume that $(M, g)$ is a smooth complete connected Riemannian manifold
and $V\in C^2(M)$ such that $\mu(d x):=e^{-V(x)} d x$ is a probability measure, where $d x$
stands for  the Riemannian volume measure, such that
  $L := \Delta - \nabla V \cdot \nabla$ generates a diffusion semigroup   $\{P_t\}_{t \geq 0}$
  on~$L^2(\mu)$, which means   that
  $\{P_t\}_{t \geq 0}$ is a symmetric strongly continuous
operator semigroup on~$L^2(\mu)$ and its generator coincides with $L$
on smooth compactly supported functions,
  $P_tf\ge 0$ if $f\ge 0$ and $P_t1=1$. Consequently, $\mu$ is $P_t$-invariant,
  i.e.,
  $$
  \int P_tf\, d\mu=\int f\, d\mu \quad \forall\, t\ge 0, \ f\in L^2(\mu).
  $$
  This triple $(M, g, \mu)$ will be called below a smooth weighted Riemannian manifold.

The curvature-dimension condition $CD(K, N)$ for the operator $L$,
where $K \in \mathbb{R}$, $N \geq 1$, is described by the Bochner-type inequality
$$
\frac{1}{2}L|\nabla f|^2 - \nabla f \cdot \nabla L f \geq
K |\nabla f|^2 +
\frac{1}{N}(Lf)^2
$$
for all smooth functions $f \in C_{0}^{\infty}(M)$.
The standard references on this topic are \cite{BE} and \cite{BakryGL}.
For example, by the classical Bochner formula
the Laplace operator on an $n$-dimensional Riemannian manifold with Ricci
curvature bounded from below by $K$ satisfies the curvature-dimension condition
$CD(K, N)$ for all $N \geq n$. However, many important diffusion operators are
intrinsically of infinite dimension, for example, for $M = \mathbb{R}^n$ the standard
Ornstein--Uhlenbeck operator $L = \Delta - x\cdot \nabla$ satisfies the condition
$CD(1, \infty)$, but does not satisfy $CD(K, N)$ with a finite number~$N$.
We recall several results from \cite{Bakry_Riesz} and \cite{Li}.
Let us define the Riesz transform $\mathcal{R}_{\varrho }$ by the formula
$$
\mathcal{R}_{\varrho } := \nabla (\varrho  - L)^{-1/2}, \ \varrho  > 0.
$$

\begin{proposition}\label{pr:riesz_transfor_estimate}
Let $(M, g, \mu)$ be a smooth weighted Riemannian
manifold
satisfying the curvature-dimension condition $CD(-\varrho , \infty)$,
$\varrho  > 0$. Then, for each $p > 1$,
 there exists a constant $C_p > 0$ depending only on $p$
such that for all $f \in C_{0}^{\infty}$ one has
$$
\|\mathcal{R}_{\varrho }f\|_{p} \leq C_p \|f\|_{p}.
$$
\end{proposition}
\begin{proof}
This is a well-known result,
 first proved by Bakry in \cite{Bakry_Riesz}, for a self-contained exposition
 and an analytical approach to this fundamental estimate we refer
 the reader to the more recent work \cite{CD}, see also \cite[Theorem~1.4]{Li}.
\end{proof}

The formulation of Theorem 1.4 in \cite{Li} also includes the case $\varrho  = 0$, but we would like to notice
that one has to be careful with the definition of $\mathcal{R}_{0}$,
since the range of $\sqrt{-L}$ on $C_0^\infty$ is not dense in $L^{2}(\mu)$.

\begin{proposition}\label{pr:meyer_inequality}
Let $(M, g, \mu)$ be a smooth weighted Riemannian
manifold
satisfying the curvature-dimension condition $CD(-\varrho , \infty)$ with some
$\varrho  \geq 0$. For every $p > 1$, there exists a constant $C_p > 0$ depending only on $p$
such that for all $f \in C_{0}^{\infty}$ one has
$$
C^{-1}_{p}\|\nabla f\|_{p}
\leq
\|\sqrt{\varrho  - L}f\|_{p}
\leq
\sqrt{\varrho } \|f\|_p + C_{q}\|\nabla f\|_p,
$$
where $1/p + 1/q = 1$ and $C_p, C_q$ are the constants provided
by Proposition~{\rm\ref{pr:riesz_transfor_estimate}}.
\end{proposition}
\begin{proof}
For $\varrho  > 0$ this is the statement of Theorem 5.5 from \cite{Li} and
in fact the proof of these inequalities was presented in \cite{Li} only in this
 case and the constants $C_p, C_q$ do not depend on $\varrho$.
The operator $-L$ is essentially self-adjoint on $C_{0}^{\infty}$ (see, e.g.,
\cite[Corollary 3.2.2]{BakryGL}) and non-negative.
Let $\{E_{\lambda}\}_{\lambda \geq 0}$
be the projection-valued measure such that
$$
-L = \int_{[0, \infty)}\lambda \,dE_{\lambda}.
$$
Let us fix $f \in C_{0}^{\infty}$.
Since $f \in D(L)$
$$
\int_{[0, \infty)}\lambda^2 \,d\langle E_{\lambda}f, E_{\lambda}f\rangle < \infty.
$$
Then by the dominated convergence theorem
\begin{multline*}
\lim_{\varrho \to 0+}\|\sqrt{-L} f  - \sqrt{\varrho - L}f\|^{2}_2 \\
= \lim_{\varrho \to 0+}
\int_{[0, \infty)} (\sqrt{\lambda} - \sqrt{\varrho + \lambda})^2
\,d\langle E_{\lambda}f, E_{\lambda}f\rangle = 0,
\end{multline*}
or, equivalently,
$\lim\limits_{\varrho  \to 0+}\sqrt{\varrho  - L}f = \sqrt{-L}f$ in $L^{2}(\mu)$.
Now the case $\varrho = 0$ easily follows by passing to the limit
as $\varrho \to 0$.
\end{proof}

The next theorem strengthens Proposition 2.2 from \cite{Led},
 established in the case $q \in [1, 2]$.
We thank the anonymous referee for pointing out that for $q>1$
it also  follows from \cite[Proposition 2.3]{ABT}.
However, we include a short proof, because in Section~5
we refer to this proof with some modification in order to cover the negative curvature case.

\begin{theorem}\label{th:p_poincare}
Let $(M, g, \mu)$ be a smooth weighted Riemannian manifold
satisfying the curvature-dimension condition $CD(0, \infty)$.
For every $q \geq 1$, there exists a positive constant $C$ depending only
on $q$ such that for any smooth function $f$ on $M$ and any $t > 0$ one has
$$
\|f - P_t f\|_q \leq C \sqrt{t}\|\nabla f\|_q.
$$
\end{theorem}
\begin{proof}
For any bounded measurable
function $h$ on~$M$ we have the following
 reverse Poincar\'e inequality (see \cite{BakryGL}):
$$
2s |\nabla P_s h|^2 \leq P_{s}(h^2) - (P_{s}h)^2.
$$
Then by Jensen's inequality for all $q' \in [2, \infty]$ we obtain
\begin{equation}\label{eq:gradient_estimate}
\|\nabla P_{s}h\|_{q'} \leq \frac{1}{\sqrt{2s}}\|h\|_{q'}.
\end{equation}
Now one can observe that
$$
\int_{M}h(f - P_t f)\,d\mu =
- \int_0^t \int_{M}hLP_{s}f\,d\mu \,ds =
\int_0^t \int_{M}\nabla P_s h \cdot \nabla f \,d\mu \,ds
$$
and by duality it is easy to see that for any $q \in [1, 2]$ one has
$$
\|f - P_{t}f \|_{q} \leq \int_0^t \frac{1}{\sqrt{2s}}\|\nabla f\|_q \,ds =
\sqrt{2t}\|\nabla f\|_q.
$$
Now let us consider $q > 2$.
For $f \in C_{0}^{\infty}$ and $t > 0$ we have (in in~$L^{2}(\mu)$)
$$
P_{t}f - f = \int_0^{\infty} K(s, t)P_{s}\sqrt{-L}f\,ds,
$$
$$
K(s, t) := \frac{1}{\sqrt{\pi}}\biggl(
\frac{\chi_{s > t}}{(s - t)^{1/2}} -
\frac{\chi_{s > 0}}{s^{1/2}}
\biggr).
$$
One can easily check that for all $t > 0$
$$
\int_0^{\infty} |K(s,t)|\,ds = \frac{4}{\sqrt{\pi}}\sqrt{t}.
$$
Taking into account Proposition \ref{pr:meyer_inequality},
we obtain
\begin{multline*}
\|f - P_t f\|_{q} \leq \int_0^{\infty}|K(s, t)|\, \|P_{s}\sqrt{-L}f\|_q\,ds \\
 \leq  \int_0^{\infty} |K(s, t)|\, \|\sqrt{-L}f\|_q\,ds \leq C(q)\sqrt{t}\|\nabla f\|_{q},
\end{multline*}
which completes the proof.
\end{proof}

\begin{remark}\label{r:gaussian_p_poincare}
\rm
For the standard Ornstein--Uhlenbeck semigroup $\{T_t\}_{t\ge0}$
given by the Mehler formula
$$
T_{t}f (x) := \int_{\mathbb{R}^d} f(e^{-t}x + \sqrt{1 - e^{-2t}} y)\,\gamma(dy)
$$
the estimate from Theorem \ref{th:p_poincare}
can be established directly with an explicit constant:
$$
\|T_{t}f - f\|_{q} \le K_{q}c_{t}\|\nabla f\|_{q},
$$
$$
K^{q}_{q} := \int_{\mathbb{R}}|x|^{q}\,\gamma(dx),\ c_{t} :=
\int_{0}^{t}\frac{e^{-s}}{\sqrt{1 - e^{-2s}}}\,ds = \arccos(e^{-t}).
$$
Indeed, for $f \in C_{0}^{\infty}(\mathbb{R}^d)$ one has
\begin{multline*}
f(e^{-t}x + \sqrt{1 - e^{-2t}}y) - f(x) = \int_0^1\frac{d}{d\tau}f(e^{-t\tau}x + \sqrt{1 - e^{-2t\tau}}y)\,d\tau \\
= t \int_0^1\nabla f(e^{-t\tau}x + \sqrt{1 - e^{-2t\tau}}y) \cdot
\biggl(-e^{-t\tau}x + \frac{e^{-2t\tau}}{\sqrt{1 - e^{-2t\tau}}}y\biggr)\,d\tau,
\end{multline*}
\begin{multline*}
\|T_{t}f - f\|^{q}_{q} \leq \int_{\mathbb{R}^d\times\mathbb{R}^d}
|f(e^{-t}x + \sqrt{1 - e^{-2t}}y) - f(x)|^{q}\gamma(dx)\gamma(dy) \\
\leq c_{t}^{q - 1}\int_{\mathbb{R}^d\times\mathbb{R}^d}\int_{[0, 1]}
\frac{te^{-t\tau}}{\sqrt{1 - e^{-2t\tau}}} \\
\times \Bigl|
\nabla f(e^{-t\tau}x + \sqrt{1 - e^{-2t\tau}}y) \cdot
(-\sqrt{1 - e^{-2t\tau}}x + e^{-t\tau}y)
\Bigr|^{q}\,d\tau \gamma(dx)\gamma(dy)\\
= c_{t}^{q - 1}\int_0^1\frac{te^{-t\tau}}{\sqrt{1 - e^{-2t\tau}}}
\int_{\mathbb{R}^d\times\mathbb{R}^d}|\nabla f(x)\cdot y|^{q}\gamma(dy)\gamma(dx)\,d\tau \\
= K_{q}^{q}c_{t}^{q}\int_{\mathbb{R}^d}|\nabla f(x)|^{q}\gamma(dx),
\end{multline*}
as announced.
\end{remark}

Recall that the Kantorovich distance
$W_p(\mu, \nu)$ of order $p \geq 1$ between two probability measures $\mu$ and $\nu$ with finite moments of order $p$ is
defined by the formula
$$
W_p(\mu,\nu)^p=\inf_{\sigma\in \Pi(\mu,\nu)} \int_{M\times M} d(x,y)^p\, \sigma(dx, dy),
$$
where inf is taken over all measures $\sigma$ from the set $\Pi(\mu,\nu)$
of Borel probability measures on $M\times M$ having projections
$\mu$ and $\nu$ onto the first and second factors, respectively; see \cite{BK12}, \cite{B18} or \cite{V}
(the case $p=1$ was considered in \cite{Kantorovich42}, \cite{Kantorovich57}).
Now let us introduce the Hopf--Lax infimum-convolutions $(Q_s)_{s > 0}$ defined
by the formula
$$
Q_{s}\varphi(x) := \inf_{y \in M}\Bigl[\varphi(y) + d^{p}(x, y)/s\Bigr], \
x \in M, \ s > 0.
$$
The dual description of the Kantorovich metric of order $p$ is given
by the equality
$$
W^{p}_{p}(\mu, \nu) = \sup_{\varphi}\biggl(
\int_{M}Q_{1}\varphi \,d\nu -
\int_{M}\varphi\,d\mu
\biggr),
$$
where the supremum is taken over all bounded continuous functions~$\varphi$, see~\cite{V}.
Alternatively, one can use the sup-convolutions defined by
$$
\widehat{Q}_{s}\varphi(x) := \sup_{y \in M}\Bigl[\varphi(y) - d^{p}(x, y)/s\Bigr], \
x \in M, \ s > 0,
$$
$$
W^{p}_{p}(\mu, \nu) = \sup_{\varphi}\biggl(
\int_{M}\varphi\,d\mu - \int_{M}\widehat{Q}_{1}\varphi \,d\nu
\biggr).
$$
A crucial role in our considerations in the next section will be played by the
 ``infinite-dimensional'' Harnack inequality that
 states that, under
 the curvature-dimension condition $CD(0, \infty)$, for any non-negative
 square-integrable function $g$ on $M$ and all $t > 0$, $x, y \in M$ one has
 one has
\begin{equation}\label{eq:harnack}
[P_{t}g(y)]^2 \leq P_{t}(g^{2})(x)e^{d^2(x, y)/2t}.
\end{equation}
For more information on Harnack inequalities of this type,
 see \cite{W-log}, \cite{WFBook}, \cite{W10}, and \cite{BGL}.

\section{Main results}

We start with establishing a weak-type bound.
The next theorem strengthens Proposition 4.2 from \cite{Led} (see also \cite{BWS16}). Set
$$\kappa(s)=   \frac{s\log s}{s\log s+1-s},\ \ s> 1.$$
It is easy to see that the function $\kappa$ is decreasing in $s\in (1,\infty)$ with
$$
\lim_{s \to 1+}\kappa(s) = \infty, \ \lim_{s \to \infty}\kappa(s) = 1.
$$

\begin{theorem}\label{th:weak_type_bound}
Let $(M, g, \mu)$ be a smooth weighted Riemannian
manifold satisfying the   condition $CD(0, \infty)$.
For every $q \geq  1$, there exists a positive constant $C$ depending only on
$q$ such that for any probability measure
$\nu = f\cdot \mu$ with a smooth density $f$
{\rm(}or $f\in W^{q,1}(\mu)${\rm)}
and any $s > 1$ one has
$$
\sup_{u\geq s}\Bigl[u^{3/2}\log^{1/2}u\Bigr]\mu(f \geq 2u)^{3/(2r)}
\leq C\kappa^{1/2}(s)\|\nabla f\|_q W_{2}(\mu, \nu),
$$
where
$r = \frac{3q}{q + 2}$.
On the other hand, the assertion does not hold if   $\log^{1/2}$ is rep
laced by an increasing positive function
$\Phi\colon [0,\infty)\to [0,\infty)$ with
$\lim\limits_{u\to\infty }\Phi(u)\log^{-1}(u)=\infty$.
\end{theorem}
\begin{proof} (a)
The proof of the first assertion is inspired by the approach of Ledoux~\cite{Led}.
For any $t > 0$ let us represent $f$ as
$$
f = (f - P_{t}f) + P_{t}f.
$$
Then for $u > 0$ we have
$$
\mu(f \geq 2u) \leq \mu\bigl(|f - P_{t}f| \geq u\bigr) +  \mu(P_{t}f \geq u).
$$
Taking into account Theorem \ref{th:p_poincare} we obtain
the estimate
$$
\mu(f \geq 2u) \leq \frac{C_{q}t^{q/2}}{u^q}\|\nabla f\|^{q}_{q} + \mu(P_{t}f \geq u).
$$
The next step is to apply
the classical entropy inequality
$$
\int_{M}hP_{t}f\,d\mu \leq
\int_{M}P_{t}f \log P_{t}f \,d\mu +
\log \int_{M}e^{h}\,d\mu.
$$
Now let us take
$$
h := \mathbbm{1}_{F}\log P_{t}f,
$$
where
$$
F = \{P_{t}f \geq u\}.
$$
Then
$$
\int_{M}\mathbbm{1}_{F}P_{t}f\log P_{t}f\,d\mu \leq
\int_{M} P_{t}f\log P_{t}f\,d\mu +
\log \int_{M}\bigl[1 + \mathbbm{1}_{F}(P_{t}f - 1)\bigr]\,d\mu ,
$$
$$
\int_{F}P_{t}f\log P_{t}f\,d\mu \leq
\int_{M} P_{t}f\log P_{t}f\,d\mu +
\log \biggl[\int_{M}1\,d\mu + \int_{F}(P_{t}f - 1)\,d\mu \biggr].
$$
Since for every $x \geq 0$ one has
$
\log (1 + x) \leq x,
$
and $\mu$ is  a probability measure, we obtain
$$
\int_{F}P_{t}f\log P_{t}f\,d\mu \leq
\int_{M} P_{t}f\log P_{t}f\,d\mu +
\int_{F}(P_{t}f - 1)\,d\mu.
$$
Note that due to \cite[Lemma 4.2]{BOGL} or \cite[Lemma 1.11]{BGL} we have
$$
\int_{M}P_{t}f\log P_{t}f\,d\mu \leq \frac{1}{4t}W^{2}_{2}(\nu, \mu).
$$
This implies the inequality
$$
\int_{F} \big(1+P_tf[\log P_tf -1]\big) \,d\mu \leq \int_{M} P_{t}f\log  P_{t}f\,d\mu
\leq \frac{1}{4t} W_2^2(\nu,\mu) \ \forall\, t>0.
$$
Combining this estimate with the definition and monotonicity of the function $\kappa$ we obtain
$$
(u\log{u})\mu(F) \leq \frac{\kappa(u)}{4t} W_2^2(\nu,\mu)
\leq \frac{\kappa(s)}{4t}  W_2^2(\nu,\mu),\  1 < s \leq u.
$$
Consequently,
$$
\mu(f \geq 2u) \leq \frac{C_{q}t^{q/2}}{u^q}\|\nabla f\|^{q}_{q} +
\frac{\kappa(s)}{4tu\log u}  W_2^2(\nu,\mu).
$$
Optimizing in $t > 0$ we obtain the desired inequality.

(b) On the other hand,  let $\Phi\colon [0,\infty)\to [0,\infty)$ be increasing such that
\begin{equation}\label{OP2}
 \lim_{u\to\infty} \Phi(u)\log^{-1} u=\infty.
\end{equation}
For any constant $C>0$, we intend to disprove the inequality
\begin{equation}\label{OP3}
\sup_{u\geq s}\Bigl[u^{3/2}\Phi(u)\Bigr]\mu(f \geq 2u)^{3/(2r)}
\leq C\kappa^{1/2}(s)\|\nabla f\|_q W_{2}(\mu, \nu),
\end{equation}
where $\nu=f\cdot\mu$,
under the condition $CD(0,\infty)$.
To this end, we take $M=\mathbb{R}$ and
$$V(x)=c+ \int_0^{|x|} d s\int_0^s h(r)d r,\ \ x\in\mathbb R,$$
where $h\in C^\infty([0,\infty))$ is nonnegative such that $h|_{[0,1/2]}= 1, h|_{[1,\infty)}=0,$
and $c \in \mathbb{R}$ is such that $\mu(dx) := \exp^{-V(x)}d x$ is a probability measure.
Then
$V\in C^\infty(\mathbb{R})$ with $V''\ge 0$ such that
the condition $CD(0,\infty)$ holds. Moreover,   there exists a constant $c_0$ such that
\begin{equation}\label{OPO}
V(x)= |x|+c_0 \quad \hbox{wnenever }  |x|\ge 1.
\end{equation}
For all $k\ge 1$, we take
$$
f_k(x)= \delta_k \{(x-k)^+\land 1\},\ \
\delta_k=\biggl(\int_{\mathbb{R}}\{(x-k)^+\land 1\} e^{-V(x)}\, d x\biggr)^{-1}.
$$
Let $\nu_k=f_k\cdot\mu$ and take $u= \frac {1}{2} e^k$. Then there exist constants $c_1,c_2>0$ such that
\begin{align*} &2c_1 u =  c_1e^k \le \delta_k\le c_2e^k = 2 c_2 u,\\
&\mu(f_k\ge 2u)\ge c_1 e^{-k}=c_1 (2u)^{-1},\\
&\|\nabla f\|_q\le c_2 e^{(1-1/q)k}=c_2 (2u)^{1-1/q},\\
&W_2(\mu,\mu_k)\le c_2 k=c_2\log(2u),\ \ u=\frac {1} {2} e^k, k\ge 1.
\end{align*}
Thus, \eqref{OP3} with $f=f_k$  implies the bound
\begin{align*}
2^{-3/2}(2u)^{1-1/q} \Phi(u) &= u^{3/2} \Phi(u) (2u)^{-(q+2)/(2q)}\\
&\le C(2u)^{1-1/q} \log(2u),\ \
u=\frac {1}{ 2} e^k, k\ge 1.
\end{align*}
Therefore,
$\liminf_{u\to\infty} \Phi(u) \log^{-1}u <\infty$
 contrary to~\eqref{OP2}.
\end{proof}

\begin{remark}
\rm
In the Gaussian case, for any probability measure of the form
$\nu = f\cdot \gamma$ with $f \in W^{q,1}(\gamma),\ q \geq 1$ one has
$$
\gamma(f \ge 2u) \le
\inf_{t > 0} \biggl[
\frac{K^{q}_{q}\arccos^{q}(e^{-t})}{u^q}\|\nabla f\|_{q}^{q} +
\frac{\kappa(s)}{2(e^{2t} - 1)u\log{u}}W_{2}^{2}(\nu, \gamma) \biggr]
$$
for all $u \geq s > 1$.
Indeed, this follows from Remark \ref{r:gaussian_p_poincare}
along the lines of the proof of Theorem \ref{th:weak_type_bound}
where one needs to take into account that in the Gaussian case we have the bound
$$
\int_{\mathbb{R}^d}T_{t}f \log T_{t}f\,d\gamma \leq \frac{W_{2}^{2}(\nu, \gamma)}{2(e^{2t} - 1)}, \ t > 0,
$$
see \cite{BGL}.
\end{remark}

The next lemma is a reinforcement of the remark made in
\cite{Led} regarding  Claim B in the proof of
Proposition 1.3 from \cite{CintiOtto}. This simple observation
will be used in the proof of Theorem \ref{theorem:main}.

\begin{lemma}\label{le:elementary}
For all $a > 1$ and $\beta \geq 0$,
any collection $\{F_k\}_{k = 1}^{\infty}$ of subsets of $M$
and any non-empty set $I\subset \mathbb{Z}_{+}$ one has
$$
\sum\limits_{k \in I}(1 + k)^{\beta}a^k \mathbbm{1}_{F_k}(x) \le
 \frac{a}{a-1} \sup_{k \in I}\big\{(1 + k)^{\beta}a^{k} \mathbbm{1}_{F_k}(x)\big\},\ x \in M.
$$
\end{lemma}
\begin{proof}
It suffices to prove this inequality
for finite $I\subset \mathbb Z_+$ such that $x\in F_k$ for some $k\in I$. In this case let us set
$k_x:=\sup I_x$, where  $I_x:= \{k\in I\colon x\in F_k\}.$
Then for any $a>1$ and $\alpha \ge 0$
\begin{multline*}
\sum\limits_{k \in I}(1 + k)^{\beta}a^k \mathbbm{1}_{F_k}(x) =
\sum\limits_{k \in I_x}(1 + k)^{\beta}a^k \leq (1+k_x)^{\beta}  \sum_{i=0}^{k_x} a^k \\
\leq (1+k_x)^\beta \cdot \frac{a^{k_x+1}-1}{a-1} \leq
(1+k_x)^\beta \cdot \frac{a^{k_x+1}}{a-1} \\
 \leq
 \frac{a}{a-1} \sup_{k\in I} \big\{(1+k)^\beta a^k \mathbbm{1}_{F_k}(x)\big\},
\end{multline*}
as announced.
\end{proof}

\begin{lemma}\label{le:kantorovich}
Let $a > 1,  \beta \geq 0$ and  $p \geq 1$. Then for any collection
 $\{F\}_{k = 1}^{\infty}$ of Borel subsets of $M$,
 every finite non-empty set $I\subset \mathbb{Z}_{+}$ and
 $\varepsilon > 0$ one has
$$
\int \varphi f\,d\mu \le \frac{1}{\varepsilon}W_{p}^{p}(\mu, f\mu)
+ \frac{a}{a-1}\int \psi_{\varepsilon}\,d\mu,
$$
where
$$
\varphi(x) = \sum_{k \in I}(1 + k)^{\beta}a^{k} \mathbbm{1}_{F_k}(x),
$$
$$
\psi_{\varepsilon}(x) = \sum_{k \in I}
\sup_{y \in M_{k, \varepsilon}(x)}(1 + k)^{\beta}a^{k}
\mathbbm{1}_{F_k}(y),
$$
$$
M_{k, \varepsilon}(x) = \Bigl\{y \in M\colon\, d(x, y)^{p} \leq
 \frac{a}{a-1}(1 + k)^{\beta}a^{k}\varepsilon \Bigr\}.
$$
\end{lemma}
\begin{proof}
By the Kantorovich duality for every $\varepsilon > 0$ we have
\begin{equation}\label{eq:k_duality}
\int \varphi f\,d\mu \leq \frac{1}{\varepsilon}W_{p}^{p}(\mu, f\mu) + \int \widehat{Q}_{\varepsilon}\varphi \,d\mu,
\end{equation}
$$
\widehat{Q}_{\varepsilon}\varphi(x) = \sup_{y \in M}\Bigl[\varphi(y) - \frac{1}{\varepsilon}d(x,y)^p\Bigr].
$$
Applying Lemma \ref{le:elementary} and letting $\Lambda:=\frac {a}{a-1}$ we obtain
\begin{multline*}
\widehat{Q}_{\varepsilon}\varphi(x) =
\sup_{y \in M}\Bigl[\sum_{k \in I}(1 + k)^{\alpha}a^{k}\mathbbm{1}_{F_k}(y)
- \frac{1}{\varepsilon}d(x,y)^p\Bigr]
\\
\leq \sup_{y \in M}\Bigl[\Lambda\sup_{k \in I}(1 + k)^{\alpha}a^{k}\mathbbm{1}_{F_k}(y)
- \frac{1}{\varepsilon}d(x,y)^p\Bigr]
\\
= \Lambda \sup_{y \in M}\sup_{k \in I}\Bigl[(1 + k)^{\alpha}a^{k}\mathbbm{1}_{F_k}(y)
- \frac{1}{\Lambda\varepsilon}d(x,y)^p\Bigr]
\\
= \Lambda \sup_{k \in I}\sup_{y \in M}\Bigl[(1 + k)^{\alpha}a^{k}\mathbbm{1}_{F_k}(y)
- \frac{1}{\Lambda\varepsilon}d(x,y)^p\Bigr]
\\
\leq \Lambda \sup_{k \in I}\sup_{y \in M_{k, \varepsilon}(x)}(1 + k)^{\alpha}a^{k}\mathbbm{1}_{F_k}(y)
\leq \Lambda \sum_{k \in I}\sup_{y \in M_{k, \varepsilon}(x)}(1 + k)^{\alpha}a^{k}\mathbbm{1}_{F_k}(y),
\end{multline*}
which together with (\ref{eq:k_duality}) completes the proof.
\end{proof}

Now we are ready to prove Theorem \ref{theorem:main}. We use Ledoux's strategy from \cite{Led} with some modifications.
For the reader's convenience all the details are presented. The key difference with the  proof
of Ledoux concerns the application of
Harnack's inequality (\ref{eq:harnack}), where the parameters are ``balanced''
in such a way that the exponential term  is no longer bounded by a constant.

\begin{proof}[Proof of Theorem~\ref{theorem:main}] (a) We first prove the claimed inequality.
For $k \in \mathbb{Z}_{+}$ and $t_k>0$
let us set
$$
A_k = \bigl\{2^k \leq f < 2^{k + 1}\bigr\},\
f_k = \min((f - 2^k)_{+}, 2^k),
$$
$$
F_k = \bigl\{P_{t_k}f_k \geq 2^{k - 1}\bigr\}, \
r = \frac{3q}{q + 2}, \
\alpha = \frac{1}{3}.
$$
For $k_0 \geq 1$  one has
\begin{multline}\label{eq:replace_integral_with_sum}
\int_{M}(f - 2^{k_0 + 1})_{+}^r\bigl(1 + \log^{r\alpha}
(1 + (f - 2^{k_0 + 1})_{+})\bigr)\,d\mu\\
=
\sum_{k = k_0 + 1}^{\infty}\int_{A_k}(f - 2^{k_0+ 1})_{+}^r\bigl(1 + \log^{r\alpha}(1 + (f - 2^{k_0 + 1})_{+})\bigr)\,d\mu \\
\leq
\sum_{k = k_0 + 1}^{\infty}2^{r(k + 1)}(1 + \log^{r\alpha}(1 + 2^{k + 1}))\mu(A_k) \\
\leq C
\sum_{k = k_0}^{\infty}(1 + k)^{r\alpha}2^{rk}\mu(f_k \geq 2^{k}),
\end{multline}
where we have used the inclusion
$
A_{k} \subseteq \{f_{k - 1} \geq 2^{k - 1}\}.
$
The constant $C$ depends only on $q$ and $k_0$.
Now let us observe that for any sufficiently smooth non-negative function $g$
(which need not be a probability density) and all $u > 0$
we have
\begin{multline}\label{eq:useful_inequality}
\mu(g \geq 2u) \leq \mu(g \geq 2u, P_{t}g \leq u) +
\mu(g \geq 2u, P_t g \geq u) \\
\leq \mu\bigl(|g - P_{t}g| \geq u\bigr) +
\frac{1}{2u}\int_{M}\mathbbm{1}_{F}g\,d\mu \\
\leq \frac{C t^{q/2}}{u^q}\|\nabla g\|^{q}_{q}
+\frac{1}{2u}\int_{M}\mathbbm{1}_{F}g\,d\mu,
\end{multline}
where $F := \{P_{t}g \geq u\}$ and the constant $C$ depends only on $q$.
Applying (\ref{eq:useful_inequality}) to $g = f_k$ and $u = 2^{k - 1}$,
$t = t_k$ we obtain the bound
\begin{multline*}
\mu(f_k \geq 2^{k}) \leq \frac{C t^{q/2}_{k}}{2^{qk}}\int_{A_k}|\nabla f|^{q}\,d\mu
+ \frac{1}{2^k}\int_{M}\mathbbm{1}_{F_k}f_{k}\,d\mu\\
\leq
\frac{C t^{q/2}_{k}}{2^{qk}}\int_{A_k}|\nabla f|^{q}\,d\mu
+ \frac{1}{2^k}\int_{M}\mathbbm{1}_{F_k}f\,d\mu,
\end{multline*}
where the inequality $f_k \leq f$ has been used.
For any $k_1 \geq k_0$ one has
\begin{multline}\label{eq:s_i_inequality}
S_{I} := \sum_{k \in I}(1 + k)^{r\alpha}2^{rk}\mu(f_{k} \geq 2^{k}) \\
 \leq
C\sum_{k \in I}(1 + k)^{r\alpha}2^{(r - q)k}t^{q/2}_{k}\int_{A_k}|\nabla f|^{q}\,d\mu +
\int_{M}\varphi f\,d\mu,
\end{multline}
where
$$I := \{k_0, \ldots, k_1\}, \
\varphi := \sum_{k \in I}(1 + k)^{r\alpha}2^{(r - 1)k}\mathbbm{1}_{F_k}.
$$
Let $\varepsilon$ be a fixed positive number.
Now we can apply Lemma \ref{le:kantorovich} to
$a = 2^{r - 1}, \ \beta = \alpha r, \ p = 2$ and obtain the inequality
\begin{equation}
\int \varphi f\,d\mu \leq \frac{1}{\varepsilon}W^{2}_{2}(f\mu, \mu) +
\frac{a}{a - 1} \int_{M} \psi_{\varepsilon}\,d\mu,
\end{equation}
where
$$
\psi_{\varepsilon}(x) =
\sum_{k \in I}\sup_{y \in M_{k, \varepsilon}(x)}(1 + k)^{r\alpha}2^{(r - 1)k}\mathbbm{1}_{F_k}(y),
$$
$$
M_{k, \varepsilon}(x) = \Bigl\{y \in M\colon d(x, y)^{2}
\leq \frac{a}{a - 1}(1 + k)^{r\alpha}2^{(r - 1)k}\varepsilon\Bigr\},
$$
Set
$$
\eta := \frac{1}{\log2} \cdot \frac{a}{a - 1}.
$$
Let us take $t_k$ such that
$$
(1 + k)^{r\alpha}2^{(r - q)k} t_{k}^{q/2} =
 \eta^{q/2}\varepsilon^{q/2},
$$
that is,
$$
t_{k} = \frac{1}{\log 2}\frac{a}{a-1}(1 + k)^{-2r\alpha/q}2^{2(1 - r/q)k}\varepsilon.
$$
By the definition of the set $F_k$ we have the pointwise estimate
$$
\mathbbm{1}_{F_k}(y) \leq
2^{-2k + 2}\bigl(P_{t_k}f_{k}(y)\bigr)^2, \ y \in M.
$$
Harnack's inequality (\ref{eq:harnack}) for the function $f_k$
and $y \in M_{k, \varepsilon}(x)$ yields
\begin{multline*}
\mathbbm{1}_{F_k}(y) \leq
2^{-2k + 2}\bigl(P_{t_k}f_{k}(y)\bigr)^2 \\
\leq
2^{-2k + 2}P_{t_k}f^{2}_{k}(x)
\exp\Bigl(
\frac{1}{2t_{k}}\frac{a}{a - 1}(1 + k)^{r\alpha}2^{(r - 1)k}\varepsilon
\Bigr).
\end{multline*}
Due to our choice of $t_k$ we have
$$
r\alpha  + 2r\alpha/q  = \alpha(r + 2r/q) =
\frac{1}{3}\Bigl(\frac{3q}{q + 2} + \frac{6}{q + 2}\Bigr) = 1,
$$
$$
(r - 1)k - 2(1 - r/q)k = (r + 2r/q - 3)k
=
\Bigl(\frac{3q}{q + 2} +  \frac{6}{q + 2} - 3\Bigr)k  = 0.
$$
This leads to the equality
$$
\exp\Bigl(
\frac{1}{2t_{k}}\frac{a}{a - 1}(1 + k)^{r\alpha}2^{(r - 1)k}\varepsilon
\Bigr) =
2^{(1 + k)/2}.
$$
Then
$$
\sup_{y \in M_{k, \varepsilon}(x)}\mathbbm{1}_{F_k}(y) \leq
2^{-2k + 2}
P_{t_k}f^{2}_{k}(x)
2^{(1 +k)/2} = 2^{-3k/2 + 5/2}P_{t_k}f^{2}_{k}(x).
$$
Since
$$
\int_{M}P_{t_k}f^{2}_{k}\,d\mu =
\int_{M} f^{2}_{k}\,d\mu \leq
2^{2k}\mu(f \geq 2^{k}) =
 2^{2k}\mu(f_{k - 1} \geq 2^{k - 1}),
$$
we have
\begin{multline*}
\int\sup_{M_{k, \varepsilon}(x)}(1 + k)^{r\alpha}2^{(r - 1)k}\mathbbm{1}_{F_k}\,d\mu \\
\leq
(1 + k)^{r\alpha}2^{(r-1)k}2^{k/2 + 5/2}\mu(f_{k - 1} \geq 2^{k - 1}) \\
=
(1 + k)^{r\alpha}2^{rk}2^{-k/2 + 5/2}\mu(f_{k - 1} \geq 2^{k - 1}).
\end{multline*}
Consequently,
\begin{multline*}
\int \varphi f\,d\mu \leq \frac{1}{\varepsilon}W^{2}_{2}(f\mu, \mu) +
\frac{a}{a - 1} \int_{M} \psi_{\varepsilon}\,d\mu \\
=
\frac{1}{\varepsilon}W^{2}_{2}(f\mu, \mu) +
\frac{a}{a - 1} \sum_{k \in I}\int_{M}\sup_{M_{k, \varepsilon}(x)}(1 + k)^{r\alpha}2^{(r - 1)k}\mathbbm{1}_{F_k}\,d\mu \\
\leq
\frac{1}{\varepsilon}W^{2}_{2}(f\mu, \mu) +
\frac{a}{a - 1}
 \sum_{k \in I}
(1 + k)^{r\alpha}2^{rk}2^{-k/2 + 5/2}\mu(f_{k - 1} \geq 2^{k - 1}).
\end{multline*}
Combining this with (\ref{eq:s_i_inequality}) and taking into account that
$$
(1 + k)^{r\alpha}2^{(r - q)k} t_{k}^{q/2} =
 \eta^{q/2}\varepsilon^{q/2},
$$
we finally have the following bound for $S_I$:
\begin{multline}\label{eq:self_improving_estimate}
S_{I} = \sum_{k \in I}(1 + k)^{r\alpha}2^{rk}\mu(f_{k} \geq 2^{k})\\
\leq C \eta^{q/2}\varepsilon^{q/2} \sum_{k \in I}\int_{A_k}|\nabla f|^{q}\,d\mu
+ \frac{1}{\varepsilon}W^{2}_{2}(f\mu, \mu)\\
+
\frac{a}{a - 1}\sum_{k \in I}(1 + k)^{r\alpha}2^{rk}2^{-k/2 + 5/2}
\mu(f_{k - 1} \geq 2^{k - 1}) \\
\leq
C \eta^{q/2}\varepsilon^{q/2} \int |\nabla f|^{q}\,d\mu
+
\frac{1}{\varepsilon}W^{2}_{2}(f\mu, \mu)\\
+ k_{0}^{r\alpha}2^{r(k_0 - 1) - k_0 + 3}
\mu(f \geq 2^{k_0}) \\+
C\frac{a}{a - 1}\sum_{k \in I}(1 + k)^{r\alpha}2^{rk}2^{-k/2 + 5/2}
\mu(f_{k} \geq 2^{k}).
\end{multline}
where the constant $C$ depends only on $q$, since
$$
(1 + 1/k)^{r\alpha} \leq 2^{r\alpha}, \ k \geq 1.
$$
We can assume that $k_0$ is sufficiently large and for all
$k \geq k_0$
$$
C\frac{a}{a - 1}2^{-k/2 + 5/2} \leq \frac{1}{2}.
$$
Then the last term on the right-hand side of (\ref{eq:self_improving_estimate}) can be replaced with~$S_I/2$:
$$
S_{I} \leq C'\varepsilon^{q/2} \int |\nabla f|^{q}\,d\mu
+ \frac{1}{\varepsilon}W^{2}_{2}(f\mu, \mu) +
C'\mu(f \geq 2^{k_0})  + \frac{1}{2}S_{I},
$$
or, equivalently,
$$
\frac{1}{2}S_{I} \leq
C' \varepsilon^{q/2} \int |\nabla f|^{q}\,d\mu
+ \frac{1}{\varepsilon}W^{2}_{2}(f\mu, \mu) +
C'\mu(f \geq 2^{k_0}),
$$
the constant $C'$ depends only on $q$ and $k_0$.
Optimizing in $\varepsilon$ and using the weak-type bound
(\ref{th:weak_type_bound}) we obtain the inequality
$$
S^{3/2r}_{I} \leq C(q, k_0) \|\nabla f\|_q W_{2}(\mu, \nu).
$$
Passing to the limit $k_1 \to \infty$ and using (\ref{eq:replace_integral_with_sum})
we derive the claimed inequality.

(b) For the second assertion, we let $\Phi$ be an increasing function such that
\begin{equation}\label{AB}
\lim_{u\to\infty} \Phi(u)\log^{-2/3}(1+u)=\infty.
\end{equation}
For any constant $C>0$, we intend to disprove the inequality
\begin{equation}\label{FI}
\bigl \|(f - C)_{+}\bigl(1 + \Phi( (f - C)_{+})\bigr)\bigr\|^{3/2}_{r}
\leq C \|\nabla f\|_{q}W_{2}(\mu, \nu),\  \ \nu=f\cdot\mu.
\end{equation}
To this end, let $M=\mathbb{R}$ and let
$V$, $\mu$,  $\nu_k := f_k\cdot\mu$
 be taken as in Step~(b) of the proof of Theorem~\ref{th:weak_type_bound}.
By \eqref{OPO},
we have $\mu= e^{V}\, d x$
and $\nu_k=f_k\cdot \mu$, hence we can
find constants $c_2>c_1>0$ such that for all $k\ge 1$ one has
$$
 c_1 e^{k}\le \delta_k\le c_2 e^k,
 \quad
W_2(\mu, \nu_k)\le c_2 k,
\quad
\|\nabla f_k\|_q\le c_2 e^{(1-1/q)k},
$$
$$
\bigl \|(f - C)_{+}\bigl(1 + F( (f - C)_{+})\bigr)
\bigr\|^{3/2}_{r}\ge c_1 e^{(1-1/q)k} \{F(e^k)\}^{3/2}.
$$
Therefore, inequality \eqref{FI} implies that
$$\limsup_{k\to\infty}k^{-1} \{\Phi(e^k)\}^{3/2} <\infty,$$
which  contradicts \eqref{AB}.
\end{proof}

\section{The Kantorovich distance $W_{1}$}

In this section we show that Theorem 1.1 from \cite{Led} admits a generalization
to the case $N = \infty$, although, unlike Theorem \ref{theorem:main} above,
where the case $p = 2$ was considered, the inequality does not include
any extra logarithmic factors. It might be possible that this result
can be further improved, we leave this question for future research.
It
would be also interesting to find a unified proof
of these inequalities covering the full scale of the Kantorovich metrics $W_p$ with $p \geq 1$.
Recall that the dual representation of the Kantorovich distance $W_{1}$
is given by the formula
\begin{equation}\label{eq:kantorovich_duality_1}
W_{1}(\mu, \nu) = \sup_{\varphi}\int_{M}\varphi \,d(\mu - \nu),
\end{equation}
where the supremum is taken over all bounded $1$-Lipschitz functions~$\varphi$, see, e.g., \cite{B18} or \cite{V}.
The next theorem is a generalization of Proposition~4.3 from \cite{Led}.

\begin{theorem}\label{th:weak_type_bound_1}
Let $(M, g, \mu)$ be a smooth weighted Riemannian
manifold satisfying the condition $CD(0, \infty)$.
For each $q \geq  1$, there exists a positive constant $C$ depending only on
$q$ such that for any probability measure
$\nu = f\cdot \mu$ with a smooth density $f$
{\rm(}or $f\in W^{q,1}(\mu)${\rm)}
one has
$$
\sup_{u > 0} \Bigl[ u^{2}\cdot\mu\bigl(|f -1| \geq 2u\bigr)^{2/r}\Bigr]
\leq C\|\nabla f\|_q W_{1}(\mu, \nu), \ r = \frac{2q}{q + 1}.
$$
\end{theorem}
\begin{proof}
For each $t > 0$ let us represent $f - 1$ as
$$
f - 1 = f -  P_{t}f + P_{t}(f - 1).
$$
Then for $u > 0$ we have
\begin{multline*}
\mu\bigl(|f - 1| \geq 2u\bigr)
\leq \mu\bigl(|f -  P_{t}f| \geq u\bigr) +
\mu\bigl(|P_t (f - 1)| \geq u\bigr) \\
\leq \mu\bigl(|f - P_{t}f| \geq u\bigr) +
\frac{1}{u}\int_{M}(\mathbbm{1}_{F_{+}} - \mathbbm{1}_{F_{-}})P_{t}(f - 1)\,d\mu
 \\
\leq \frac{C t^{q/2}}{u^q}\|\nabla f\|^{q}_{q}
+ \frac{1}{u}\int_{M}(\mathbbm{1}_{F_{+}} - \mathbbm{1}_{F_{-}})P_{t}(f - 1)\,d\mu \\
=
\frac{C t^{q/2}}{u^q}\|\nabla f\|^{q}_{q}
+\frac{1}{u}\int_{M}P_{t}(\mathbbm{1}_{F_{+}} -\mathbbm{1}_{F_{-}}) \,d(\nu - \mu),
\end{multline*}
where
$$
F_{+} :=\bigl\{P_{t}(f - 1) \geq u\bigr\},
\ F_{-} := \bigl\{P_{t}(f - 1) \leq -u\bigr\}
$$
and $C$ depends only on $q$.
By the gradient estimate (\ref{eq:gradient_estimate}) we have
$$
|\nabla P_{t}(\mathbbm{1}_{F_{+}} - \mathbbm{1}_{F_{-}})|
\leq \frac{1}{\sqrt{2t}},
$$
hence the Kantorovich duality (\ref{eq:kantorovich_duality_1})
yields the bound
$$
\int_{M}P_{t}(\mathbbm{1}_{F_{+}} -\mathbbm{1}_{F_{-}}) \,d(\nu - \mu)
\leq \frac{1}{\sqrt{2t}}\frac{1}{u}W_{1}(\mu, \nu).
$$
Finally, we have
$$
\mu\bigl(|f - 1| \geq 2u\bigr) \leq \frac{C t^{q/2}}{u^q}\|\nabla f\|^{q}_{q}
+ \frac{1}{\sqrt{2t}}\frac{1}{u}W_{1}(\mu, \nu),
$$
so optimizing in $t > 0$ we arrive at the desired inequality.
\end{proof}

From the weak-type bounds provided by Theorem \ref{th:weak_type_bound_1} we deduce the corresponding strong ones.

\begin{theorem}\label{th:strong_type_bound_1}
Let $(M, g, \mu)$ be a smooth weighted Riemannian
manifold satisfying the  condition $CD(0, \infty)$.
For each $q >  1$, there exists a positive constant $C$ depending only on
$q$ such that for any probability measure $\nu = f\cdot \mu$ with a smooth density $f$
{\rm(}or $f\in W^{q,1}(\mu)${\rm)}
one has
$$
\bigl \|(f - C)_{+}\bigr\|^{2}_{r}
\leq C \|\nabla f\|_{q}W_{1}(\mu, \nu), \ r = \frac{2q}{q + 1}.
$$
\end{theorem}
\begin{proof}
Following the lines of the proof of
Theorem \ref{theorem:main}, we introduce
$$
A_k = \bigl\{2^k \leq f < 2^{k + 1}\bigr\},\
f_k = \min((f - 2^k)_{+}, 2^k),
$$
$$
F_k = \bigl\{P_{t_k}f_k \geq 2^{k - 1}\bigr\}, \
r = \frac{2q}{q + 1}, \ a := 2^{r - 1}
$$
and establish the inequality
$$
\int_{M}\bigl(f - 2^{k_0 + 1}\bigr)_{+}^r\,d\mu
\leq C
\sum^{\infty}_{k = k_0}2^{rk}\mu(f_k \geq 2^{k}),
$$
where $C$ depends only on $q$ and $k_0$.
Next we bound $\mu(f_k \geq 2^{k})$ using the estimate
$$
\mu(f_k \geq 2^{k}) \leq \frac{C t^{q/2}_{k}}{2^{qk}}\int_{A_k}|\nabla f|^{q}\,d\mu
+ \frac{1}{2^k}\int_{M}\mathbbm{1}_{F_k}f\,d\mu.
$$
Then for
$$I := \{k_0, \ldots, k_1\}, \
\varphi := \sum_{k \in I}2^{(r - 1)k}\mathbbm{1}_{F_k}.
$$
we obtain the chain of inequalities
\begin{multline*}
S_{I} := \sum_{k \in I}2^{rk}\mu(f_{k} \geq 2^{k}) \\
 \leq
C\sum_{k \in I}2^{(r - q)k}t^{q/2}_{k}\int_{A_k}|\nabla f|^{q}\,d\mu +
\int_{M}\varphi f\,d\mu \\
\leq
C\sum_{k \in I}2^{(r - q)k}t^{q/2}_{k}\int_{A_k}|\nabla f|^{q}\,d\mu + \frac{1}{\varepsilon}W_{1}(\mu, \nu) \\+
\frac{a}{a - 1}\sum_{k \in I}\int \sup_{M_{k, \varepsilon}(x)}2^{(r - 1)k}
\mathbbm{1}_{F_k}\,d\mu,
\end{multline*}
where
$$
M_{k, \varepsilon}(x) := \Bigl\{y \in M\colon d(x, y) \leq \frac{a}{a - 1}2^{(r -1)k}\varepsilon\Bigr\}.
$$
Let us set
$$
\eta := \frac{1}{\log2} \cdot \frac{a}{a - 1}
$$
and take $t_k$ such that
$$
2^{(r - q)k} t_{k}^{q/2} = \eta^{q/2}\varepsilon^{q/2},
$$
that is,
$$
t_{k} = \frac{1}{\log 2}\frac{a}{a-1}2^{2(1 - r/q)k}\varepsilon.
$$
Then
$$
\mathbbm{1}_{F_k}(y) \leq
2^{-2k + 2}\bigl(P_{t_k}f_{k}(y)\bigr)^2, \ y \in M,
$$
and by Harnack's inequality applied to the function $f_k$
and $y\in M_{k, \varepsilon}(x)$ we have
\begin{multline*}
\mathbbm{1}_{F_k}(y) \leq
2^{-2k + 2}\bigl(P_{t_k}f_{k}(y)\bigr)^2 \\
\leq
2^{-2k + 2}P_{t_k}f^{2}_{k}(x)
\exp\Bigl(
\frac{1}{2t_{k}}\frac{a}{a - 1}2^{(r - 1)k}\varepsilon
\Bigr).
\end{multline*}
The definition of $t_k$ yields the equality
$$
(r - 1)k - 2(1 - r/q)k = (r + 2r/q - 3)k
=
\Bigl(\frac{3q}{q + 2} +  \frac{6}{q + 2} - 3\Bigr)k  = 0.
$$
This leads to the equality
$$
\exp\Bigl(
\frac{1}{2t_{k}}\frac{a}{a - 1}2^{(r - 1)k}\varepsilon
\Bigr) =
2^{1/2}.
$$
Then
$$
\sup_{y \in M_{k, \varepsilon}(x)}\mathbbm{1}_{F_k}(y) \leq
2^{-2k + 2}
P_{t_k}f^{2}_{k}(x)
2^{1/2} = 2^{-2k + 5/2}P_{t_k}f^{2}_{k}(x).
$$
Recalling that
$$
\int_{M}P_{t_k}f^{2}_{k}\,d\mu =
\int_{M} f^{2}_{k}\,d\mu \leq
2^{2k}\mu(f \geq 2^{k}) =
 2^{2k}\mu(f_{k - 1} \geq 2^{k - 1}),
$$
similarly to the proof of Theorem \ref{theorem:main}
we obtain the ``recursive'' inequality for all sufficiently large $k_0$:
$$
S_{I} \leq C \eta^{q/2}\varepsilon^{q/2}\int |\nabla f|^q\,d\mu +
\frac{1}{\varepsilon}W_{1}(\mu, \nu) +
2^{r(k_0 - 1) - k_0 + 3}\mu\bigl(f \geq 2^{k_0}\bigr) + \frac{1}{2}S_I.
$$
Taking into account Theorem \ref{th:weak_type_bound_1}
 it is easy to complete the proof.
\end{proof}

Now let us consider the case $q =1$.

\begin{theorem}\label{th:strong_type_bound_1_1}
Let $(M, g, \mu)$ be a smooth weighted Riemannian
manifold satisfying the  condition $CD(0, \infty)$.
Then for each probability measure $\nu = f\cdot \mu$ with a smooth density $f$
{\rm(}or $f\in W^{1,1}(\mu)${\rm)}
one has
$$
\bigl \|f -1\bigr\|^{2}_{1}
\leq 2 \|\nabla f\|_{1}W_{1}(\mu, \nu).
$$
\end{theorem}
\begin{proof}
This inequality is a particular case of the results from \cite{BWS15}.
We present a proof here for the reader's convenience.
For a smooth function $g \in C_{0}^{\infty}$ with $\|g\|_{\infty} \leq 1$ we have
\begin{multline*}
\biggl|\int_{M} (f - 1)g\,d\mu \biggr|
= \biggl|\int_{M} (f - 1)P_{t}g\,d\mu - \int_{M} (f - 1) \int_{[0, t]}\frac{d}{ds}P_{s}g\,d\mu\biggr|
\\
\leq W_{1}(\nu, \mu)\|\nabla P_tg\|_{\infty} + \int_{0}^{t}\int |\nabla f|\, |\nabla P_{s}g|\,d\mu\,ds
\\
\leq W_{1}(\nu, \mu)\|\nabla P_{t}g\|_{\infty}
+ \|\nabla f\|_{1} \int_{[0, t]}\|\nabla P_{s}g\|_{\infty}\,ds
\\
\leq
W_{1}(\nu, \mu) \frac{1}{\sqrt{2t}} +
 \|\nabla f\|_{1} \int_{[0, t]}\frac{1}{\sqrt{2s}}\,ds =
 W_{1}(\nu, \mu) \frac{1}{\sqrt{2t}}  + \sqrt{2t}\|\nabla f\|_1.
\end{multline*}
Optimizing in $t > 0$ we obtain the desired inequality.
\end{proof}

Of course, this theorem covers the case of the standard Gaussian measure $\gamma_d$ on $\mathbb{R}^d$ with the usual metric.
As already mentioned in our note \cite{BWS15}, in this case the following two-sided inequality holds
for all functions $f$ from the Gaussian Sobolev class $W^{1,1}(\gamma_d)$ having zero integral against~$\gamma_d$:
$$
\frac{\|f\|_{L^1(\gamma_d)}^2}{2\|\nabla f\|_{L^1(\gamma_d)}} \le \|f\cdot\gamma_d\|_K \le \|\nabla f\|_{L^1(\gamma_d)}.
$$

In \cite[Proposition 1 and Proposition 2]{BWS15}
we constructed two examples showing that the bound from this theorem can fail with any constant
if $\mu$ does not satisfy the indicated condition (in one of these examples $\mu$
is a measure on the real line with the usual
metric and in the other example $M$ is a two-dimensional complete connected Riemannian submanifold
in~$\mathbb{R}^d$
and $\mu$ is its Riemannian volume).

The obtained inequalities involving the
Kantorovich distance $W_1$ of order~$1$
can be combined with the estimate provided by Theorem 1.1 from~\cite{Led}.

\begin{proposition}
Let $(M, g, \mu)$ be a smooth Riemannian
manifold satisfying the curvature-dimension condition $CD(0, N)$.
Then, given $p,q\ge 1$, there exists a constant $C>0$ depending only on $p,q,N$ such that
for each probability measure $\nu = f\cdot \mu$ with a smooth density $f$
{\rm(}or $f\in W^{1,1}(\mu)${\rm)}
one has
$$
\|f-1\|_r^\theta \le C\Big(\|\nabla f\|_q W_p(\mu,\nu) +
\bigl[\|\nabla f\|_1W_1(\mu,\nu)\bigr]^{\theta/(2r)}\Bigr),
$$
where
$$
r= \frac{1 + \frac{1}{p} + \frac{1}{N}}{\frac{1}{p} + \frac{1}{q}}, \
\theta= 1+\frac{1}{p} +\frac{1}{N}.
$$
In particular, for $p=q=1$ this becomes
\begin{equation}\label{AA}
\|f-1\|_{1+\frac{1}{2N}}^{2+\frac{1}{N}}
\le C\|\nabla f\|_1W_1(\mu, \nu).
\end{equation}
\end{proposition}
\begin{proof}
By \cite[Theorem 1.1]{Led}
there exists a constant $C > 1$ such that
\begin{equation}\label{eq:ledoux_inequality_N}
 \|(f-C)_{+}\|_{r}^\theta \le C \|\nabla f\|_q W_p(\mu,\nu)
\end{equation}
and by Theorem \ref{th:strong_type_bound_1_1}
we have
\begin{equation}\label{eq:our_inequality_1_1}
\|f-1\|_1 \leq \sqrt{2}\|\nabla f\|_1 W_1(\mu,\nu).
\end{equation}
Since $f$ is non-negative,
we have the following trivial bound:
$$
|f - 1|^r \leq 2^{r - 1}\bigl[(f - C)_{+}\bigr]^{r} + (2C)^{r - 1}|f - 1|.
$$
Integrating this inequality over $M$ and applying (\ref{eq:ledoux_inequality_N}) and (\ref{eq:our_inequality_1_1})
we get the desired bound.
\end{proof}

\begin{remark}
\rm
 The inequality \eqref{AA}
is sharp in the sense that, whenever $C>0$ and  $\Phi\colon [0,\infty)\to [0,\infty)$
is an increasing function  with
\begin{equation} \label{BB} \lim_{u\to\infty} \Phi(u)u^{-2-\frac 1 N} =\infty,
\end{equation}
the following inequality does not hold:
\begin{equation}\label{AA1}
\Phi(\|f-1\|_{1+\frac{1}{2N}})
\le C\|\nabla f\|_1W_1(\mu, \nu).
\end{equation}
Indeed, let $M=(\mathbb S^1)^N$, where $\mathbb S^1$ is the unit circle,
which is equivalent to $[0,2\pi)$ with the periodic boundary.
For every $n\ge 1$ we set $h_n(s)= \min\{ns, (2-ns)^+\}, s\in [0,2\pi)$ and
$$f_n(x)= \prod_{i=1}^N h_n(x_i),\ \ x=(x_1,\cdots, x_N)\in [0,2\pi)^N.$$
Let $\nu_n = f_n\cdot d x$,
 $\mu= \frac{1} {2\pi}\, d x$,
We have $W_1(\nu_n,\mu)\le 2\pi$,
and  there exist   constants $c_2\ge c_1>0$   such that
$$\|f_n-1\|_{1+\frac{1} {2N}}\ge c_1 n^{\frac {N} {2N+1}},
\ \ \|\nabla f_n\|_1\le c_2 n,\ \ n\ge 1.$$
Thus, \eqref{AA1} implies that
$\liminf_{u\to\infty} \Phi(u)u^{-2+\frac 1 N} <\infty$,
which contradicts \eqref{BB}.
\end{remark}

\section{Extensions to the negative curvature case }

In this section we briefly discuss some extensions for negatively
curved weighted Riemannian manifolds.
We assume that $(M, g, \mu)$
satisfies the curvature-dimension condition $CD(-\varrho , \infty)$,
$\varrho  > 0$ and that additionally the logarithmic Sobolev inequality holds:
\begin{equation}\label{ineq:log_sobolev}
\int_{M}f^2 \log f^2 \,d\mu \leq \frac{2}{\lambda}
\int_{M}|\nabla f|^2\,d\mu
\end{equation}
for all
 $ f \in C^{1}(M)$ with
 $$
 \int_{M}f^2\,d\mu = 1.
$$
For example, according to  \cite{W04}, under the curvature-dimension condition
$CD(-\varrho, \infty)$ the finiteness of the integral
$$
\int_{M}\exp\bigl(\varepsilon d^2(x_0, x)\bigr)\,d\mu
$$
for some $x_0 \in M$ and $\varepsilon > \varrho/2$
ensures the validity of
the log-Sobolev inequality \ref{ineq:log_sobolev}.
The main idea of the considerations below is that even though in this case the curvature bound alone
does not guarantee the required semigroup estimates,
nevertheless they can be established under some additional assumptions about $(M, g, \mu)$.

\begin{proposition}\label{pr:n_curv_gradient_bound}
Let $(M, g, \mu)$ be a smooth weighted Riemannian
satisfying the
curvature-dimension condition $CD(-\varrho, \infty)$ with some
$\varrho \geq 0$. Assume also that
the log-Sobolev inequality {\rm(\ref{ineq:log_sobolev})} holds for some
$\lambda > 0$. Then for each $p \in [2, \infty)$ there exists $C > 0$
depending only on $\varrho, \lambda, p$ such that
$$
\|\nabla P_{t}h\|_{p} \leq \frac{C}{\sqrt{t}}\|h\|_{p},
\ t > 0, \ h \in L^{p}(M, \mu).
$$
\end{proposition}
\begin{proof}
Using the standard approximation arguments one can see
that it is sufficient to establish this inequality just for
$h \in C_{b}(M)$.
Applying \cite[Corollary~4.2]{W-log}
we obtain the inequality
$$
|\nabla P_{t}h| \leq \frac{C}{\sqrt{t \wedge 1}}
\bigl(
P_{t}|h|^p
\bigr)^{\frac{1}{p}}.
$$
This implies the bound
\begin{equation}\label{eq:neg_curv_p}
\|\nabla P_{t}h\|_{p} \leq \frac{C}{\sqrt{t \wedge 1}}\|h\|_p.
\end{equation}
Now one can observe that the log-Sobolev inequality (\ref{ineq:log_sobolev})
ensures that the generator of the semigroup $\{P_t\}_{t\ge0}$
has a spectral gap larger or equal to~$\lambda$, in particular,
for any $\varphi \in L^2(\mu)$ with zero integral against $\mu$ one has
$$
\|P_{t}\varphi\|_{2} \leq e^{-\lambda t}\|\varphi\|_2.
$$
Consequently, for $h \in C_{b}(M)$ we have
\begin{equation}\label{eq:neg_curv_p_2}
\|\nabla P_{t}h\|_{2} \leq \frac{Ce^{-\lambda t}}{\sqrt{t \wedge 1}}\|h\|_2.
\end{equation}
Combining estimates (\ref{eq:neg_curv_p}), (\ref{eq:neg_curv_p_2})  and
applying the standard interpolation theorem we obtain the desired inequality.
\end{proof}

\begin{theorem}\label{pr:n_curv_poincare}
Let $(M, g, \mu)$ be a smooth weighted Riemannian
manifold satisfying  the
curvature-dimension condition $CD(-\varrho, \infty)$ with some
$\varrho \geq 0$. Assume also that
the log-Sobolev inequality {\rm(\ref{ineq:log_sobolev})} holds for some
$\lambda > 0$. Then, for each $q \in (1, 2]$, there exists $C > 0$
depending only on $\varrho, \lambda, p$ such that
for all $f\in W^{q,1}(\mu)$ one has
$$
\|f - P_{t}f\|_{q} \leq C \sqrt{t}\|\nabla f\|_q,\ t > 0.
$$
\end{theorem}
\begin{proof}
This follows from Proposition \ref{pr:n_curv_gradient_bound}
along the lines of the proof of Theorem \ref{th:p_poincare}.
\end{proof}

\begin{theorem}\label{pr:n_curv_entropy_bound}
Let $(M, g, \mu)$ be a smooth weighted Riemannian
manifold satisfying the
curvature-dimension condition $CD(-\varrho, \infty)$ with some
$\varrho \geq 0$. Assume also that
the log-Sobolev inequality {\rm(\ref{ineq:log_sobolev})} holds for some
$\lambda > 0$. Then  there exists $C > 0$
depending only on $\varrho, \lambda, p$ such that
for every probability density $f$ one has
$$
\int_{M}P_{t}f \log P_{t}f\,d\mu \leq \frac{C}{t}W^2_{2}(\mu, f\cdot\mu).
$$
\end{theorem}
\begin{proof}
According to \cite{W10} the curvature condition $CD(-\varrho, \infty)$ implies the log-Harnack
inequality
$$
P_{t}\log g (x) \leq \log P_{t}g(y) +
\frac{\varrho \cdot d^2(x,y)}{2(1 - e^{-2\varrho t})}.
$$
Applying this inequality to $g := P_{t}f$ and
integrating  with respect to the optimal coupling
of the measures $f\cdot \mu$  and $\mu$
(see, e.g., \cite[Corollary~1.2]{RW10}) we get the bound
$$
\int_{M}P_{t}f \log P_{t}f\,d\mu \leq \frac{\varrho}{2(1 - e^{-2\varrho t})}
W^{2}_{2}(\mu, f\cdot \mu).
$$
Consequently, for all $t \in (0, 1)$ we have
$$
\int_{M}P_{t}f \log P_{t}f\,d\mu \leq \frac{C}{t}W^{2}_{2}(\mu, f\cdot \mu).
$$
Next, it is known that the log-Sobolev inequality (\ref{ineq:log_sobolev})
implies the following bound for any probability density $g$ with respect to~$\mu$:
$$
\int_{M}P_{t}g\log P_{t}g\,d\mu \leq e^{-\lambda t}\int_{M}g\log g\,d\mu.
$$
Applying this inequality to $g = P_{1}f$, we obtain the estimate
\begin{multline*}
\int_{M}P_{t}f\log P_{t}f\,d\mu
\leq e^{-\lambda (t - 1)}\int_{M}P_{1}f\log P_{1}f\,d\mu \\
\leq C e^{-\lambda t}W^{2}_{2}(\mu, f\cdot\mu), \  t \geq 1.
\end{multline*}
Now it is easy to complete the proof.
\end{proof}

\begin{theorem}
Let $(M, g, \mu)$ be a smooth weighted Riemannian
manifold  satisfying the
curvature-dimension condition $CD(-\varrho, \infty)$ with some
$\varrho \geq 0$. Assume also that
the log-Sobolev inequality {\rm(\ref{ineq:log_sobolev})} holds for some
$\lambda > 0$. Then for each $q \in (1, 2]$ there exists $C > 0$
depending only on $\varrho, \lambda, q$ such that
for every smooth probability density $f$
{\rm(}or $f\in W^{q,1}(\mu)${\rm)}
and every $s > 1$ one has
$$
\sup_{u\geq s}\Bigl[u^{3/2}\log^{1/2}u\Bigr]\mu(f \geq 2u)^{3/(2r)}
\leq C\kappa^{1/2}(s)\|\nabla f\|_q W_{2}(\mu, \nu),
$$
where
$$
r = \frac{3q}{q + 2}, \
\kappa(s) := \frac{s\log s}{s\log s + 1 - s}, \ s > 1.
$$
\end{theorem}
\begin{proof}
This follows from  Theorem \ref{pr:n_curv_poincare} and Theorem \ref{pr:n_curv_entropy_bound}
along the lines of the proof of Theorem \ref{th:weak_type_bound}.
\end{proof}

Let us conclude this section with a generalization of Theorem \ref{th:strong_type_bound_1_1}.
For a function $f\in L^1(\mu)$ let
$$\|f-\mu(f)\|_K=\sup_{g\in C^\infty(M), \ \|\nabla g\|_\infty\le 1} \int_M fg \, d\mu.$$
When $\nu=f\cdot \mu$ is a probability measure, we have $W_1(\mu,\nu)=\|f-1\|_K$.

\begin{theorem}
Let $(M, g, \mu)$ be a smooth weighted Riemannian
manifold satisfying the
curvature-dimension condition $CD(-\varrho, \infty)$,
$\varrho \geq 0$. Assume also that
the semigroup $\{P_t\}_{t\ge0}$ satisfies the inequality
\begin{equation}\label{eq:strong_ergodicity}
\|P_{t}g\|_{\infty} \leq ce^{-\lambda t} \|g\|_{\infty}, \
 \int_{M}g\,d\mu = 0, \ t \geq 0.
\end{equation}
with some $c, \lambda > 0$.
Then there exists $C > 0$ depending only on $\varrho, c, \lambda$
such that for every integrable smooth function $f$ with zero integral against $\mu$ one has
$$
\|f\|^{2}_{1} \leq C\|\nabla f\|_{1} \|f\|_{K}.
$$
\end{theorem}
\begin{proof}
First, let us remind that for any $h \in C_{b}(M)$ we have the
pointwise inequality (see \cite[Corollary 4.2]{W-log})
$$
|\nabla P_{t}h| \leq \frac{C}{\sqrt{t \wedge 1}}\bigl(P_{t}|h|^2\bigr)^{1/p},
$$
consequently,
$$
\|\nabla P_{t}h\|_{\infty} \leq \frac{C}{\sqrt{t \wedge 1}}\|h\|_{\infty}.
$$
Next, using our additional assumption about the semigroup $\{P_t\}_{t\ge0}$ it readily seen
that for $t \geq 1$ one has
$$
\|\nabla P_{t}h\|_{\infty} = \|\nabla P_{1}P_{t - 1}h\|_{\infty}\leq
C'e^{-\lambda t}\|h\|_{\infty}.
$$
Combining these two bounds we obtain
$$
\|\nabla P_{t}h\|_{\infty} \leq \frac{C}{\sqrt{t}}\|h\|_{\infty}
$$
and,  consequently,
$$
\|f - P_{t}f\|_{1} \leq C\sqrt{t}\|\nabla f\|_1.
$$
Now it is easy to complete the proof similarly to
Theorem \ref{th:strong_type_bound_1_1}, see also our short note~\cite{BWS15}.
\end{proof}

\begin{remark}
\rm
According to \cite{W04}, the log-Sobolev inequality and the strong ergodicity
(inequality (\ref{eq:strong_ergodicity})) are incomparable, but both follow from the ultraboundedness:
$\|P_t\|_{1\to\infty}<\infty$ for $t>0$.
See also \cite{RW03} for more details.
\end{remark}

We thank two anonymous referees for useful comments.

This research was supported  by
the Russian Science Foundation grant  17-11-01058 (Sections 2 and~3)
and the NNSFC grants 11771326, 11831014, and 11921001 (Sections~4 and~5).
{\sloppy

}

\end{document}